\documentclass[11pt]{article}
\usepackage{graphicx}
\usepackage{pgfplots}
\usepackage[dvipsnames]{xcolor}
\usepackage[colorlinks=true,citecolor=blue]{hyperref}
\usepackage{natbib}
\usepackage{graphicx}
\usepackage{amsfonts}
\usepackage{amsmath}
\usepackage{amssymb}
\usepackage{url}
\usepackage{fancyhdr}
\usepackage{indentfirst}
\usepackage{enumerate}
\usepackage{titlesec}
\usepackage{amsthm}
\usepackage{dsfont}
\usepackage[misc]{ifsym}
\usepackage{enumitem}

\theoremstyle{definition}

\newtheorem{example}{Example}

\theoremstyle{plain}
\newtheorem{theorem}{Theorem}
\newtheorem{lemma}{Lemma}
\newtheorem{proposition}{Proposition}
\newtheorem{corollary}{Corollary}
\theoremstyle{remark}

\usepackage{cases}

\theoremstyle{definition}

\def\E{\mathbb{E}}

\def\R{\mathbb{R}}

\definecolor{proofblue}{HTML}{236B91}
\definecolor{prooforange}{HTML}{C56E26}

\usepackage[onehalfspacing]{setspace}

\usepackage{bm}
\usepackage{tikz-qtree}
\usepackage{tikz}

\def\id{\mathds{1}}

\pgfplotsset{compat=1.18}

\title{Martingale representations without finite or well-defined mean}

\author{Benjamin C\^ot\'e\thanks{Department of Statistics and Actuarial Science, University of Waterloo, Canada.   \texttt{b3cote@uwaterloo.ca}} \and  Lilian Hu\thanks{Department of Statistics and Actuarial Science, University of Waterloo, Canada.   \texttt{l66hu@uwaterloo.ca}} \and  Ruodu Wang\thanks{Department of Statistics and Actuarial Science, University of Waterloo, Canada.   \texttt{wang@uwaterloo.ca}}}

\date{\today}

\begin{document}
	\maketitle
	\begin{abstract}
The classical Strassen theorem characterizes convex order through martingale representations under the assumption of finite first moments. We show that a natural reformulation of the theorem remains true on the real line for distributions without finite or even well-defined means. We further investigate the higher-dimensional setting, where this equivalence fails in general. Martingale representations are nonetheless guaranteed in higher dimensions when every nonzero linear projection of the second marginal has undefined mean.
\medskip

	\end{abstract}

 \noindent\rule{\textwidth}{0.5pt}

\section{Introduction}

The celebrated 
Strassen theorem \citep[Theorem~8]{S65} gives a characterization of the convex order by stating that two distributions are convexly ordered if and only if they admit a martingale representation. This equivalence is notably central to the construction of martingales \citep{HPRY11} and to the theory of martingale optimal transport \citep{BHP13}.
  The original formulation of Strassen's theorem  assumes that the random variables have finite means. This assumption de facto makes the conditional expectations from the martingale representation well defined. 
Accordingly, the convex order is often defined on $L^1$, for instance, in the textbook \cite{MS02}, purposefully to avoid the complications arising from infinite means.

Recently, \cite{CW26} showed that most of the defining properties of the convex order continue to hold when extending the order to $L^0$, the set of all random variables (with possibly infinite or undefined means). This begs the question whether the characterization of the convex order given by Strassen's theorem also continues to hold on $L^0$.

Under its usual formulation, this may not be the case because the involved conditional expectations may be undefined. 
We present a slightly different formulation of Strassen's theorem in which the conditioning $\sigma$-algebra need not be $\sigma(X)$. This formulation is equivalent to the classical one on $L^1$, and we show that it remains valid on $L^0$. The proof builds on the existence of shadow measures established by \cite[Lemma~4.6]{BJ16} and their potential representation in \cite[Theorem~4.7]{BHN22}.
Therefore, this new unified  version of 
 Strassen's theorem  (equivalent to the classic one on $L^1$) should be regarded as the canonical form of the theorem. 
The convex order considered here is the same as that in \cite{CW26}, including for random variables beyond $L^1$. For example, when a random variable has undefined mean, it dominates all others with respect to the convex order, and a Strassen representation is available for any pair of distributions where at least one has undefined mean.  

We also present a multivariate extension to the convex order on $L^0$. Under the multivariate setting, the equivalence between convex ordering and the existence of a Strassen representation
 fails to hold. We give a counterexample. We furthermore discuss the case where the mean of the dominating distribution is infinite along every ray and show that, analogously to the univariate case with undefined mean, a Strassen representation is then possible with any other distribution.   

The paper is divided as follows. In Section~\ref{sect:preliminaries}, we give some preliminaries on the convex order. 
Sections~\ref{sect:univariate} and~\ref{sect:multivariate} contain our main results on the one-dimensional and multi-dimensional Strassen theorems, respectively.

\section{Preliminaries}
\label{sect:preliminaries}

Let $(\Omega, \mathcal{F}, \mathbb{P})$ be an atomless probability space.  Define $L^0$ as the set of all random variables, 
$L^1$ as the set of random variables with finite first moment, and $L^{\infty}$ as the set of all bounded random variables.  For random variables, relations such as ``$X \leq Y$'' are interpreted to hold almost surely. Almost surely equal random variables are treated as identical.

 For any $X,Y\in L^0$, we say that $Y$ dominates $X$ in convex order, denoted by $X\leq_{\rm cx} Y$, when 
          $ \mathbb{E}[\varphi(X)] \leq \mathbb{E}[\varphi(Y)]$
    holds for all convex functions $\varphi:\mathbb{R}\to\mathbb{R}$ such that both expectations are well-defined (but possibly infinite). An expectation $\mathbb{E}[X]$ is undefined when $\mathbb{E}[X_+] = \mathbb{E}[X_-] = \infty$. 
    Note that the literature contains two versions of $\leq_{\rm cx}$ depending on whether infinite expectations are to be considered, and although both coincide on $L^1$, they differ on $L^0$ \citep[Proposition~1]{CW26}. The results in this paper exclusively concern the order taking infinite expectations into consideration, and this is aligned with the accrued coherence this conveys to the convex order, as argued in \cite{CW26}. 

Let 
$\mathcal{M}_d$ be the set of all finite measures and $\mathcal{P}_d$ be the set of all probability measures on $(\mathbb{R}^d, \mathcal{B}(\mathbb{R}^d))$, where $\mathcal{B}$ indicates the corresponding $\sigma$-algebra, $d\in\mathbb{N}=\{1,2,\ldots\}$. For each, we omit the index $d$ when $d=1$. Write $\delta_k$ for the Dirac measure at $k\in\mathbb{R}$. 
    
 We alternatively write $\mu \leq_{\rm cx} \nu$ for probability measures $\mu,\nu\in\mathcal{P}$ when $\int \varphi \,\mathrm{d}\mu \leq \int \varphi \,\mathrm{d}\nu$ holds for all convex functions $\varphi$ such that both integrals are well-defined; this is consistent with the definition for random variables. We extend this definition to general non-probability measures $\lambda,\eta\in\mathcal{M}$ by additionally requiring $\lambda(\mathbb{R}) = \eta(\mathbb{R})$. If we rather require $\lambda(\mathbb{R}) \leq \eta(\mathbb{R})$ and comparisons to hold for nonnegative convex functions only, this is the extended convex order, denoted by $\leq_{\rm ecx}$. 

A random variable $X$ with undefined mean will produce undefined or positive-infinite expectations $\mathbb{E}[\varphi(X)]$ for any nonconstant convex function $\varphi$, and therefore it dominates all random variables with respect to the convex order. We record this fact in a lemma for future referencing. 
\begin{lemma}
\label{lem:always}
   If $Y\in L^0$ has undefined mean, then $X\leq_{\rm cx} Y$ always holds, for all $X\in
    L^0$.
\end{lemma}
We will often employ the standard Cauchy distribution as an example of distribution with undefined mean. The standard Cauchy has density $f$, with respect to the Lebesgue measure, given by $f(x) = 1/(\pi(1+x^2))$ with $x\in\mathbb{R}$. 

To formulate our results, we need conditional expectations  expressed in terms of $\sigma$-algebras for random variables possibly outside $L^1$, and this may need some explanation. We use an equivalent definition to the one in \citet[p.213 and p.218]{Shiryaev1996}.
Write $x\wedge y=\min\{x,y\}$ for $x,y\in \R$, $x_{-}=\max\{-x,0\}$ and $x_+=\max\{x,0\}$. For a nonnegative random variable $Y$ and a $\sigma$-algebra $\mathcal{G}$,
we define 
$$
\mathbb{E}[Y\mid\mathcal{G}]=
\lim_{m\to\infty} \E[Y\wedge m\mid \mathcal G],
$$
where   the limit exists   by the monotone convergence theorem. 
The conditional expectation $\mathbb{E}[Y\mid\mathcal{G}]$ is a random variable taking values in $[0,\infty]$, and 
it is the a.s.~unique random variable $Z$ that is measurable with respect to $\mathcal{G}$  satisfying 
 $
\mathbb{E}\left[Z \id_A \right] = \mathbb{E}\left[Y\id_A\right]$  for every  $A\in\mathcal{G}$; see  \citet[p.218]{Shiryaev1996}.
 Clearly, if $Y$   is $\mathcal{G}$-measurable, then $\mathbb{E}[Y\mid \mathcal{G}] = Y$. For a general $Y\in L^0$, we define $$\mathbb{E}[Y\mid\mathcal G]=\mathbb{E}[Y_+\mid\mathcal G]-\mathbb{E}[Y_{-}\mid\mathcal G],$$ whenever the difference is well-defined, that is, $\min\{\mathbb{E}[Y_+\mid\mathcal G],\mathbb{E}[Y_{-}\mid\mathcal G]\}<\infty$ a.s. If  both terms are infinite on a non-null set, then $\mathbb{E}[Y\mid \mathcal{G}]$ is undefined 
 (this is similar to the usual expectation $\E[X]$, which is undefined when $\E[X_+]=\E[X_-]=\infty$). 
When we write $\mathbb{E}[Y\mid\mathcal{G}]=X$ for a real-valued random variable $X$, we implicitly require both $\mathbb{E}[Y_+\mid\mathcal{G}]$ and $\mathbb{E}[Y_-\mid\mathcal{G}]$ to be finite a.s.
 
\section{Strassen's theorem without finite mean on the real line}
\label{sect:univariate}

We first present the classical Strassen theorem \citep[Theorem~8]{S65}, and amend it slightly by allowing for more flexibility in the conditioning $\sigma$-algebra.


\begin{theorem} 
\label{th:S1}
Consider $\mu, \nu \in \mathcal{P}$ both with finite first moment. The following are equivalent:
\begin{enumerate}[label={\rm (\roman*)}, ref={\rm (\roman*)}]
    \item \label{item:S1-cx} $\mu\leq_{\rm cx} \nu$;
    \item \label{item:S1-G} there exist $X\sim \mu$, $Y\sim \nu$, and a $\sigma$-algebra $\mathcal{G}$ such that $\mathbb{E}[Y \mid \mathcal{G}] = X$.
      \item \label{item:S1-sX} there exist $X\sim \mu$, $Y\sim \nu$ such that $\mathbb{E}[Y \mid \sigma(X) ] = X$;
\end{enumerate}
\end{theorem}
\begin{proof}
The equivalence between \ref{item:S1-cx} and \ref{item:S1-sX} is the classical result \cite[Theorem 8]{S65}, which requires means to be finite. The implication \ref{item:S1-sX}$\implies$\ref{item:S1-G} is direct. For the implication \ref{item:S1-G}$\implies$\ref{item:S1-sX}: Suppose that $\mathbb{E}[Y\mid \mathcal{G}] = X$ for some $\mathcal{G}\subseteq \mathcal{F}$. This indicates that $X$ is measurable with respect to $\mathcal{G}$. Therefore, 
\begin{equation*}
    \mathbb{E}[Y\mid \sigma(X)] = \mathbb{E}\left[\mathbb{E}[Y\mid \mathcal{G}] \mid \sigma(X) \right] = \mathbb{E}[X\mid\sigma(X)] = X,  
\end{equation*}
where the first equality is given by the tower property. Note that $Y\in L^1$ means that each conditional expectation above is finite and well-defined, so that applying the tower property poses no issue. 
\end{proof}

Throughout, we will say that $\mu,\nu$ admit a \textit{Strassen representation} if \ref{item:S1-G} holds. 
In \ref{item:S1-G}, the random variable $X$ is measurable with respect to $\mathcal{G}$ and therefore $\sigma(X)\subseteq \mathcal{G}$, meaning that more information than the one provided by $X$ is required. The key point of the amendment is that the equivalence with the convex order is maintained even when one allows more information to be fed to $X$. This is precisely the observation enabling to extend Strassen's theorem to $L^0$: the information provided by $X$ alone may be so little that $\mathbb{E}[Y\mid \sigma(X)]$ may be infinite or undefined with positive probability, meaning it is not a proper random variable, and thus it automatically prevents a Strassen representation using $\sigma(X)$. We illustrate this caveat in the following example.

\begin{example}
\label{ex:Cauchy-1}
    Let $\mu$ be the Dirac measure at $0$, and $\nu$ be the standard Cauchy distribution. Because $\nu$ has undefined mean, it holds that $\mu\leq_{\rm cx} \nu$  by Lemma~\ref{lem:always}. A Strassen representation with minimal information, as in Theorem~\ref{th:S1}\ref{item:S1-sX}, would mean that there are $X=0$ and $Y\sim \nu$ such that $\mathbb{E}[Y \mid\sigma(X)] = 0$ almost surely. However, by the  degeneracy of $X$, we have
    $\mathbb{E}[Y\mid\sigma(X)] = \mathbb{E}[Y]$, which is undefined (and thus not equal to 0). Nonetheless, $\mu$ and $\nu$ admit a Strassen representation: namely, choose $\mathcal{G} = \sigma(|Y|)$. Then, for any $\omega\in\Omega$, letting $y = |Y(\omega)|$,  
    \begin{equation*}
    \mathbb{E}[Y \mid \mathcal{G} ](\omega) = \frac{yf_\nu(y) - yf_\nu(-y)}{f_{|Y|}(y)} = 0 ,\quad \text{for all }\omega\in \Omega 
    \end{equation*}
    where $f_\nu(y) = f_\nu(-y)$ holds by symmetry of the standard Cauchy distribution around $0$. Thus, $\mathbb{E}[Y \mid \mathcal{G} ]$ is degenerate at $0$, and we obtain $\mathbb{E}[Y\mid \mathcal{G}] = X$. 
\end{example}

This leads to our main result, stating that a Strassen representation is always available whenever random variables are convexly ordered, even beyond $L^1$. 

\begin{theorem} 
\label{th:S0}
Consider $\mu, \nu \in \mathcal{P}$. The following are equivalent:
\begin{enumerate}[label={\rm (\roman*)}, ref={\rm (\roman*)}]
    \item \label{item:S0-cx} $\mu \leq_{\rm cx} \nu$;
    \item \label{item:S0-G} there exist $X\sim \mu$, $Y\sim \nu$, and a $\sigma$-algebra $\mathcal{G}$ such that $\mathbb{E}[Y \mid \mathcal{G}] = X$;
    \item \label{item:S0-alphabeta} There are sequences of finite measures $\{\alpha_n:n\in\mathbb{N}\}$ and $\{\beta_n:n\in\mathbb{N} \}$ each with finite first moment such that 
    \begin{equation}
        \mu = \sum_{n=1}^{\infty} \alpha_n, \quad \nu = \sum_{n=1}^{\infty} \beta_n, \quad \text{and} \quad \alpha_n \leq_{\rm cx} \beta_n \quad \text{for every }n\in\mathbb{N}. 
    \end{equation}
\end{enumerate}
In particular, in \ref{item:S0-G}, $\mathbb{E}[Y \mid \mathcal{G}]$ is well-defined and finite almost surely. If $\nu$ has well-defined first moment, possibly infinite, one can choose $\mathcal{G} = \sigma(X)$.   
\end{theorem}
\begin{proof}
    
For the implication \ref{item:S0-cx}$\implies$\ref{item:S0-alphabeta}: By assumption $\mu \leq_{\rm cx} \nu$.
Set $\mu_0=\mu$ and $\nu_0=\nu$. For $n\geq 1$, suppose that $\mu_{n-1}\leq_{\rm cx}\nu_{n-1}$ with same nonzero total mass. Choose a compactly supported submeasure $\alpha_n\leq\mu_{n-1}$ such that $\alpha_n(\mathbb R)=\mu_{n-1}(\mathbb R)/2$. This is doable since $\mu_{n-1}\left([-M,M]\right)\rightarrow\mu_{n-1}(\R)=r_{n-1}>0$ as $M\rightarrow \infty$, we can choose $k>0$ such that $\mu_{n-1}\left([-k,k]\right)\geq r_{n-1}/2$. Hence, we can define $\alpha_n=c\mu_{n-1}\vert_{[-k,k]}$ where $c=r_{n-1}/(2\mu_{n-1}([-k,k]))\leq1$. By Lemma~\ref{lem:decompo} (below), there exists an integrable submeasure $\beta_n\leq\nu_{n-1}$ such that $\alpha_n\leq_{\rm cx}\beta_n$ and $\mu_{n-1}-\alpha_n\leq_{\rm cx}\nu_{n-1}-\beta_n$. Define $\mu_n=\mu_{n-1}-\alpha_n$ and $\nu_n=\nu_{n-1}-\beta_n$. Hence, $\mu_n\leq_{\rm cx}\nu_n$, and we can do the same procedures recursively.

Moreover, $\mu_n(\mathbb R)=\nu_n(\mathbb R)=2^{-n}\mu(\mathbb R)\to 0$ as $n\to\infty$. Since $\mu=\sum_{i=1}^n\alpha_i+\mu_n$ and $\nu=\sum_{i=1}^n\beta_i+\nu_n$ for every $n$, letting $n\to\infty$ yields $\mu=\sum_{n=1}^{\infty}\alpha_n$ and $\nu=\sum_{n=1}^{\infty}\beta_n$, as desired.

For the implication \ref{item:S0-alphabeta}$\implies$\ref{item:S0-G}: For each $n\in\mathbb{N}$, we have $\alpha_n \leq_{\rm cx} \beta_n$ by assumption, and $\beta_n$ has finite mean. Thus, the classical Strassen theorem (Theorem~\ref{th:S1}) yields that 
\begin{equation}
    \mathbb{E}[B_n \mid \sigma(A_n)] = A_n, \label{eq:strassen-alphabeta} 
\end{equation}
for some $A_n\sim \alpha_n/\alpha_n(\mathbb{R})$ and $B_n\sim \beta_n/\beta_n(\mathbb{R})$, for every $n\in\mathbb{N}$. Let $N$ denote a random variable independent of all $(A_n,B_n)$, with $\mathbb{P}(N=n) = \alpha_n(\mathbb{R}) = \beta_n(\mathbb{R})$. 
(These probability masses sum to $1$ because $\mu$ and $\nu$ are probability measures.) Note that, by assumption, there are $X\sim\mu$ and $Y\sim\nu$ that are mixtures of $A_n$ and $B_n$ respectively, where $N$ is the mixing random variables; indeed,
\begin{equation*}
 \nu = \sum_{n=1}^{\infty} \beta_n = \sum_{n=1}^{\infty} \mathbb{P}(N=n) \frac{\beta_n}{\beta_n(\mathbb{R})}  \quad \text{and} \quad
\mu = \sum_{n=1}^{\infty} \alpha_n = \sum_{n=1}^{\infty} \mathbb{P}(N=n) \frac{\alpha_n}{\alpha_n(\mathbb{R})}  . 
\end{equation*}
This yields the following stochastic representations of $Y$ and $X$: 
\begin{equation}
    Y = \sum_{n=1}^{\infty} \id_{\{N=n\}} B_n \quad \text{and} \quad X = \sum_{n=1}^{\infty} \id_{\{N=n\}} A_n. \label{eq:storep-alphabeta}
\end{equation}
Therefore, 
\begin{align}
    \mathbb{E}[Y\mid \sigma(X,N)] &= \mathbb{E}[Y_+\mid \sigma(X,N)] - \mathbb{E}[Y_-\mid \sigma(X,N)]\notag\\ &= \mathbb{E}\left.\left[\left(\sum_{n=1}^{\infty} \id_{\{N=n\}} B_n\right)_+ \;\right| \sigma(X,N) \right] - \mathbb{E}\left.\left[\left(\sum_{n=1}^{\infty} \id_{\{N=n\}} B_n \right)_-\;\right| \sigma(X,N) \right]\notag\\
    &= \sum_{n=1}^{\infty} \id_{\{N= n\}} \left(\mathbb{E}[(B_n)_+  \mid \sigma(X,N)] - 
    \mathbb{E}[(B_n)_-  \mid \sigma(X,N)]\right)\notag\\
    &=\sum_{n=1}^{\infty} \id_{\{N= n\}} \mathbb{E}[B_n  \mid \sigma(X,N)] = \sum_{n=1}^{\infty} \id_{\{N= n\}} \mathbb{E}[B_n \mid \sigma(A_n,N)], \label{eq:condexp-alphabeta}
\end{align}
where the last equality is because $X = A_n$ on $\{N=n\}$.
Furthermore note that
\begin{equation}
    \mathbb{E}[B_n \mid \sigma(A_n,N)] = \mathbb{E}[B_n \mid \sigma(A_n)] = A_n, \quad \text{for every }n\in\mathbb{N}, 
    \label{eq:indep-alphabeta}
\end{equation}
where the first equality is by independence of $N$ and $(A_n,B_n)$, and the second is from \eqref{eq:strassen-alphabeta}. Combining \eqref{eq:storep-alphabeta}, \eqref{eq:condexp-alphabeta} and \eqref{eq:indep-alphabeta} produces
\begin{equation*}
    \mathbb{E}[Y\mid \mathcal{G}] = \sum_{n=1}^{\infty} \id_{\{N=n\}}\mathbb{E}[B_n\mid\sigma(A_n,N)] = \sum_{n=1}^{\infty} \id_{N=n} A_n = X,
\end{equation*}
which yields the desired claim by choosing $\mathcal{G} = \sigma(X,N)$. 

For the implication \ref{item:S0-G}$\implies$\ref{item:S0-cx}: For any convex function $\varphi : \mathbb{R}\to \mathbb{R}$, 
\begin{equation*}
    \mathbb{E}[\varphi(X)] = \mathbb{E}\left[ \varphi\left( \mathbb{E}[Y\mid \mathcal{G}] \right)\right]  \leq \mathbb{E}\left[  \mathbb{E}[\varphi(Y)\mid \mathcal{G}]\right] = \mathbb{E}[\varphi(Y)],
\end{equation*}
provided both $\mathbb{E}[\varphi(X)]$ and $\mathbb{E}[\varphi(Y)]$ are well-defined, where the inequality is Jensen's inequality. This is $X\leq_{\rm cx} Y$ by definition. 

Lastly, we want to show that in \ref{item:S0-G}, if $\nu$ has well-defined first moment, possibly infinite, one can choose $\mathcal{G} = \sigma(X)$. We consider the case where $\E[Y^-]<\infty$, and the case $\E[Y^+]<\infty$ is similar. From $\E[Y\mid G]=X$, we obtain that $\E[Y^+\mid G]+X^-=X^++\E[Y^-\mid G]$. From the tower property of the nonnegative random variables, we have that $\E[Y^+\mid \sigma(X)]+X^-=X^++\E[Y^-\mid \sigma(X)]$. Hence $\E[Y^+\mid \sigma(X)]=X+\E[Y^-\mid \sigma(X)]$. Since $\E[Y^-\mid \sigma(X)]<\infty$ a.s. from $\E[Y^-]<\infty$, we obtain that $\E[Y\mid \sigma(X)]=X$
\end{proof}

Theorem~\ref{th:S0}
suggests, in particular, that the equivalence 
\ref{item:S0-cx}$\Longleftrightarrow$\ref{item:S0-G}
should be viewed as the canonical form of Strassen's theorem, because it works both on $L^0$ and $L^1$, and it immediately gives the classic form  on $L^1$.  Moreover, in the standard definition of a two-step martingale $(X,Y)$, the $\sigma$-algebra $\mathcal G$ such that $\E[Y\mid \mathcal G]=X$ should be general and need not be $\sigma(X)$, so this version of Strassen's theorem is naturally aligned with the canonical formulation   of martingales.


The proof of Theorem~\ref{th:S0} relies on the following Lemma~\ref{lem:decompo}, which allows us to extract an integrable pair of measures satisfying convex order while preserving convex order for the remaining measures. 



\begin{lemma}
\label{lem:decompo}
    Consider $\lambda, \eta \in \mathcal{M}$ with same nonzero total mass, and any $\alpha\leq\lambda$ with compact support and $\alpha(\mathbb{R})<\lambda(\mathbb{R})$. If $\lambda\leq_{\rm cx}\eta$, then there exists a measure $\beta \leq \eta$ with well-defined finite first moment and satisfying
    \begin{equation*}
        \alpha\leq_{\rm cx} \beta \quad \text{and} \quad \lambda - \alpha \leq_{\rm cx} \eta - \beta. 
    \end{equation*}
\end{lemma}
\begin{proof}
    The proof is provided in Appendix~\ref{proof:decompo}. 
\end{proof}
The main idea of the proof is to rely on previous results on the existence and construction 
of shadow measures. If $\eta$ has finite mean, we directly choose
 $\beta$ as the shadow measure of $\alpha$ in $\eta$, which gives $\alpha\leq_{\rm cx}\beta\leq\eta$, and $\beta$'s existence is guaranteed by \citet[Lemma~4.6]{BJ16}. When $\eta$ has infinite or undefined mean, we first truncate it to have finite mean and choose $\beta$ to be the shadow measure under this truncated one. The remaining convex order $\lambda - \alpha \leq_{\rm cx} \eta - \beta$ is then established using some properties of the shadow measures, such as its representation in \citet [Theorem~4.7]{BHN22} and its minimality, together with some properties of infinite-mean or undefined-mean distributions.
The noteworthiness of Lemma~\ref{lem:decompo} lies in the finite first moment of the measure $\beta$, which then enables to apply the classical version of Strassen's theorem to this portion of the measure.

The general theory of \cite{S65} requires random variables to form a separable Banach space. The assumption of finite mean in the classical Strassen theorem  directly comes from this requirement, since there is no norm $||\cdot||$ on $L^0$ compatible with convergence in probability that would make $(L^0, ||\cdot||)$ a Banach space. Existing extensions of \cite{S65}'s theory keep this requirement, sometimes indirectly. For example,  \citet[Theorem 2.6.1]{MS02} for $\leq_{\rm cx}$, exclude infinite-mean and undefined-mean random variables, as the integral-order approach explicitly removes them from the ordered space.

Combining Theorem~\ref{th:S0} and Lemma~\ref{lem:always} directly yields the following probabilistic result. 
\begin{corollary}
\label{th:WOW}
    Consider any $\nu\in \mathcal{P}$ with undefined mean. For every $\mu\in \mathcal{P}$, there exists a $\sigma$-algebra $\mathcal{G}\subseteq \mathcal{F}$ such that 
   $
   \mathbb{E}[Y\mid \mathcal{G}] = X
    $
    for some $Y\sim\nu$ and $X\sim \mu$.
\end{corollary}

The existence result above does not give an explicit construction of the $\sigma$-algebra $\mathcal{G}$. 
The following proposition presents a construction for the case where $X$ is degenerate and $Y$ is absolutely continuous. 

\begin{proposition}
\label{prop:WOW-constant}
    Let $k\in\mathbb{R}$, and consider any $\nu \in\mathcal{P}$ that has undefined mean and density $f_{\nu}$. Then, if $Y\sim\nu$, it holds that $\mathbb{E}[Y\mid\sigma(\lfloor h_k(Y)\rfloor)] = k$, where $h_k(y) = \int_k^y (x-k) f_{\nu}(x)\,\mathrm{d}x$. 
\end{proposition}

\begin{proof}
    
The function $h_k$ is nonnegative. It is  increasing on $y>k$ and decreasing on $y<k$ because $f_{\nu}$ is nonnegative; moreover, $\{h_k(y) : y\leq k \} = \{ h_k(y) : y \geq k\} = [0,\infty)$. One notices that $h_k$ corresponds to portions of the mean of $Y$ above or below $k$, and thus this last observation holds only because $Y$ has undefined mean.
As a result,  the following values are well-defined for each $n \in\mathbb{N}_0$:
\begin{equation*}
a_n = \sup \{ y \leq k: h_k(y) = n \} \quad \text{and} \quad b_n =  \inf \{y\geq k : h_k(y) = n\}. 
\end{equation*}
Define events $A_n = \{ n \leq h_k(Y) < n+1 \}$ for each $n\in\mathbb{N}$. Because of the shape of $h_k$, 
\begin{equation*}
    A_n = \{ \omega\in\Omega : Y(\omega) \in (a_{n+1}, a_n]\cup [b_n, b_{n+1})\}, \quad n\in\mathbb{N}_0. 
\end{equation*}
Also,
\begin{align*}
    \mathbb{E}[(Y-k)\id_{A_n}] &= \int_{a_{n+1}}^{a_n} (x-k) f_{\nu}(x)\,\mathrm{d}x + \int_{b_{n}}^{b_{n+1}} (x-k) f_{\nu}(x)\,\mathrm{d}x\\ &= h_k(a_{n}) - h_k(a_{n+1}) + h_k(b_{n+1}) - h_k(b_n)  = n - (n+1) + (n+1) - n = 0. 
\end{align*}
Hence, $\mathbb{E}[Y\id_{A_n}] = k \mathbb{P}(A_n) $ for every $n\in\mathbb{N}_0$.
Note that $\sigma(\lfloor h_k(Y)\rfloor) = \sigma(\{ A_n : n\in\mathbb{N}_0\})$ and $\{ A_n : n\in\mathbb{N}_0\}$ forms a partition of the sample space with $\mathbb P(A_n)>0$ for every $n\in\mathbb N_0$. Therefore,
\begin{equation*}
    \mathbb{E}[Y \mid \sigma(\lfloor h_k(Y)\rfloor)]  = \mathbb{E}[Y\mid \sigma(\{ A_n : n\in\mathbb{N}\})] =\sum_{n=0}^\infty \id_{A_n}\frac{\mathbb{E}[Y\id_{A_n}]}{\mathbb{P}(A_n)} = k,  
\end{equation*}
as desired.
\end{proof}

A direct observation from Proposition~\ref{prop:WOW-constant} is that the tower rule $\mathbb{E}\left[\mathbb{E}[Y \mid \mathcal{G}] \right] = \mathbb{E}[Y]$ is invalid in general for $Y \in L^0$: the left-hand side may be well-defined while the right-hand side is not, as in Example~\ref{ex:Cauchy-1}. Furthermore, by taking $k_1, k_2\in\mathbb{R}$,  $k_1 \neq k_2$, Proposition~\ref{prop:WOW-constant} yields
\begin{equation*}
    \mathbb{E}[Y\mid \sigma(\lfloor h_{k_1}(Y) \rfloor) ] \neq  \mathbb{E}[Y \mid \sigma(\lfloor h_{k_2}(Y) \rfloor)].
\end{equation*}
This is impossible if $Y$ were integrable, since the tower property would require $k_1=k_2=\mathbb E[Y]$.  




\section{Strassen's theorem in higher dimensions}
\label{sect:multivariate}

A multivariate version of the convex order is as follows: 
    For any $X,Y\in (L^0)^d$, the random vector $Y$ dominates $X$ in convex order, denoted by $X\leq_{\rm cx} Y$, when 
    $\mathbb{E}[\varphi(X)] \leq \mathbb{E}[\varphi(Y)]$
    holds for all convex functions $\varphi: \mathbb{R}^d \to \mathbb{R}$ such that both expectations are well-defined (but possibly infinite). In terms of measures, for $\mu,\nu\in\mathcal{M}_d$, we write $\mu \leq_{\rm cx} \nu$ when $\int \varphi\, \mathrm{d}\mu \leq \int \varphi\, \mathrm{d}\nu$ holds for all convex $\varphi: \mathbb{R}^d \to \mathbb{R}$ whenever both integrals are well-defined and $\mu(\mathbb{R}^d) = \nu(\mathbb{R}^d$).
Note that infinite expectations need to be considered in the comparisons above, and thus this definition of the multivariate convex order is the direct extension to the order discussed so far. 

\begin{proposition}
\label{th:S1-multi}
    In Theorem~\ref{th:S1}, the equivalences continue to hold when letting $\mu,\nu\in \mathcal{P}_d$ and assuming that both measures have finite first moment in each dimension.  
\end{proposition}
\begin{proof}
The statement of \citet[Theorem~8]{S65} already allows for $\mu, \nu \in \mathcal{P}_d$, and assumes that both have finite means in each dimension. The argument in the proof of Theorem~\ref{th:S1} for the legitimacy of feeding more information to Strassen representations remains valid on $(L^1)^d$. 
\end{proof}

A similar statement cannot be made without the assumption of finite first moments: convexly ordered random vectors may not admit a Strassen representation. The following example illustrates this. 

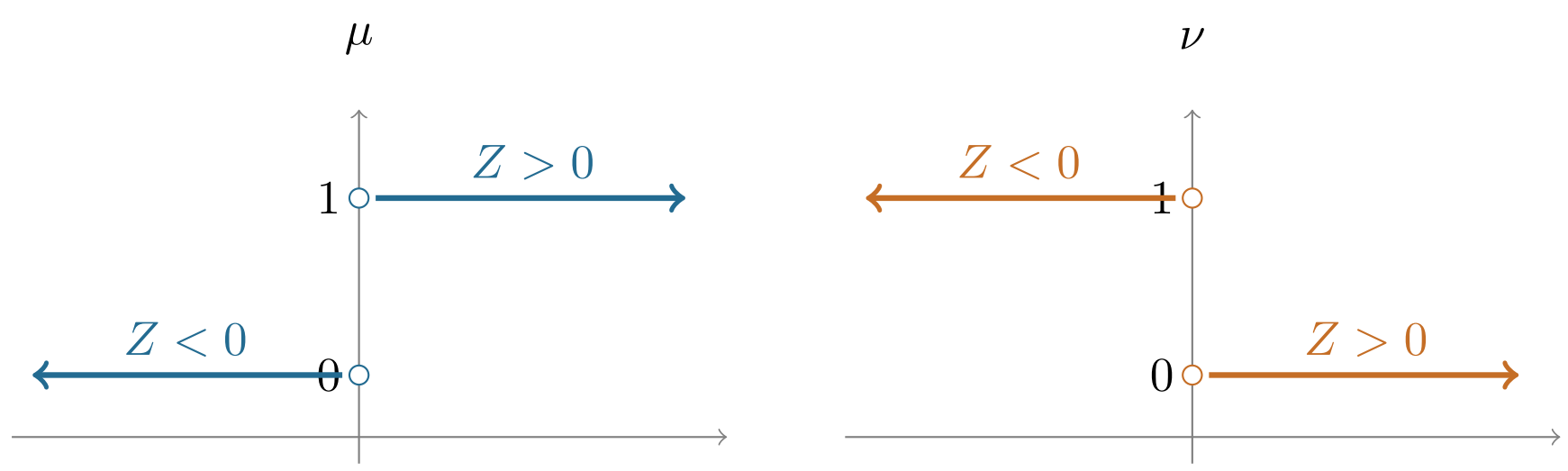
\begin{figure}[tb]
\centering
\begin{tikzpicture}[x=1.02cm,y=1.3cm,font=\small]
 \begin{scope}
  \node at (0,1.9){\textbf{$\mu$}};
  \draw[->,gray] (-2.5,-.35)--(2.65,-.35) ;
  \draw[->,gray] (0,-.5)--(0,1.5) ;
  \node[left] at (0,0){$0$}; \node[left] at (0,1){$1$};
  \draw[very thick,proofblue,->] (.12,1)--(2.35,1);
  \draw[very thick,proofblue,->] (-.12,0)--(-2.35,0);
  \draw[proofblue,fill=white] (0,1) circle(2pt);
  \draw[proofblue,fill=white] (0,0) circle(2pt);
  \node[proofblue,above] at (1.25,1){$Z>0$};
  \node[proofblue,above] at (-1.25,0){$Z<0$};
 \end{scope}
 \begin{scope}[shift={(6,0)}]
  \node at (0,1.9){\textbf{$\nu$}};
  \draw[->,gray] (-2.5,-.35)--(2.65,-.35) ;
  \draw[->,gray] (0,-.5)--(0,1.5) ;
  \node[left] at (0,0){$0$}; \node[left] at (0,1){$1$};
  \draw[very thick,prooforange,->] (-.12,1)--(-2.35,1);
  \draw[very thick,prooforange,->] (.12,0)--(2.35,0);
  \draw[prooforange,fill=white] (0,1) circle(2pt);
  \draw[prooforange,fill=white] (0,0) circle(2pt);
  \node[prooforange,above] at (-1.25,1){$Z<0$};
  \node[prooforange,above] at (1.25,0){$Z>0$};
 \end{scope}
 \end{tikzpicture}
\caption{Support of the distributions from Example~\ref{ex:no-strassen-multi}. Distributions $\mu$ and $\nu$ convexly dominate each other, but they do not admit a Strassen representation.}
\label{fig:no-strassen}
\end{figure}

\begin{example}
\label{ex:no-strassen-multi}
    Let $\mu$ be the distribution of $(Z,\id_{\{Z>0\}})$ and $\nu$ be the distribution of $(Z,\id_{\{Z < 0\}})$, where $Z$ follows a standard Cauchy distribution. A representation of the supports of $\mu$ and $\nu$ is provided in Figure~\ref{fig:no-strassen}. First, let us argue that $\mu\leq_{\rm cx}\nu$ and $\nu\leq_{\rm cx} \mu$. For any convex function $\varphi:\mathbb{R}^2\to\mathbb{R}$, one of two possibilities can occur: either
    \begin{enumerate}[nosep]
        \item  $\varphi(z_1,z_2)=g(z_2)$ for some convex $g$. Then $\int \varphi \,\mathrm{d}\mu = \int \varphi \,\mathrm{d}\nu$; or,
        \item $\varphi(z_1,z_2) \geq az_1 + b z_2 + c$ for some $a\neq 0$, for all $z_1,z_2\in\mathbb{R}$. If $a>0$, then $\int\varphi_+\,\mathrm{d}\mu \geq \mathbb{E}[(aZ + b\id_{\{Z>0\}} + c)_+] =  \infty$, and similarly $\int\varphi_+\,\mathrm{d}\nu = \infty$. If $a<0$, the infinite negative first moment of the Cauchy distribution also implies $\int\varphi^+\,\mathrm d\mu=\int\varphi^+\,\mathrm d\nu=+\infty$. Hence, whenever $\int\varphi\,\mathrm d\mu$ and $\int\varphi\,\mathrm d\nu$ are well-defined, they must both equal $+\infty$.
    \end{enumerate}
    Therefore, $\int \varphi \,\mathrm{d}\mu = \int \varphi\, \mathrm{d}\nu$ always holds whenever the integrals are well-defined, so $\mu\leq_{\rm cx} \nu$ and $\nu \leq_{\rm cx} \mu$. However, for any $X=(X_1,X_2)\sim\mu$ and $Y=(Y_1,Y_2)\sim\nu$, it is impossible to write $\mathbb{E}[Y\mid\mathcal{G}] = X$ for any $\sigma$-algebra $\mathcal{G}$. Suppose this is the case, then this would mean $\mathbb{E}[Y_2\mid \mathcal{G}] = X_2$. 
    Since $\mathbb E[Y_2\mid\mathcal G]=X_2$ and both $X_2$ and $Y_2$ take values in ${0,1}$, we first show that $Y_2=X_2$ almost surely. Since $X_2$ is $\mathcal G$-measurable, $\mathbb E[Y_2 \id_{A_0}]=\mathbb E[\mathbb E[Y_2\mid\mathcal G]\id_{A_0}]=\mathbb E[X_2 \id_{A_0}]=0$. Since $Y_2 \id_{A_0}\geq0$, we obtain $Y_2=0$ a.s.~on $A_0$. When $X_2=1$, apply the similar arguments on $1-Y_2$, we obtain that $Y_2=1$ a.s.~on $\{X_2=1\}$. Therefore, $Y_2=X_2$ a.s.

Recall that $X_2=\id_{\{X_1>0\}}$ and $Y_2=\id_{\{Y_1<0\}}$. Thus, on $A=\{X_1>0\}$, we have that $Y_1<0$. Since $X_1$ is $\mathcal G$-measurable, $\mathbb E[Y_1 \id_A\mid\mathcal G]=\id_A\mathbb E[Y_1\mid\mathcal G]\leq0$. However, the assumed martingale condition implies $\mathbb E[Y_1 1_A\mid\mathcal G]=1_A X_1=X_1>0$ on $A$. This is a contradiction. 
\end{example}



Without finite means, the disconnection between convex order and martingalization arising in dimensions $d\geq 2$  has been previously remarked by \cite{EH92}, which they discuss through fusions of distributions; however, their version of the convex order removes infinite expectations. Therefore, in their Example~3.19, the random vectors $X$ and $Y$ are not convexly ordered under our definition,
which means that their example does not provide a counterexample in our setting. Example~\ref{ex:no-strassen-multi} above works under both settings. 


Although convex order alone does not guarantee a Strassen representation in higher dimensions, we now identify a condition on $\nu$ that ensures such a representation for every $\mu$. We say that a distribution on $\mathbb{R}^d$ has \emph{infinite mean along every ray} when 
\begin{equation}
    \int (w^{\mathrm{T}}y)_+ \nu(\mathrm{d}y) = \infty \quad \text{for every $w\in\mathbb{R}^{d}\backslash \{0\}$}.
    \label{eq:all-directions}
\end{equation}
Equivalently, every nonzero linear projection of $\nu$ has undefined mean. 
For example, the joint distribution of independent standard Cauchy random variables has infinite mean along every ray. 
This property is the analogue to having undefined mean on $\mathbb{R}$ in that it characterizes umbrella-distributions dominating all other distributions in $\mathbb{R}^d$ in terms of the convex order. 

\begin{theorem}
\label{th:umbrella-multi}
    For any $\nu\in\mathcal{P}_d$, the following are equivalent: 
    \begin{enumerate}[label={\rm (\roman*)}, ref={\rm (\roman*)}]
        \item \label{item:multi-infinitemean} $\nu$ has infinite mean along every ray;
        \item \label{item:multi-cx} $\mu \leq_{\rm cx} \nu$ holds for every $\mu\in\mathcal{P}_d$; 
        \item \label{item:multi-G} for every $\mu\in\mathcal{P}_d$, there exist $X\sim\mu$, $Y\sim \nu$ and a $\sigma$-algebra $\mathcal{G}$ such that $\mathbb{E}[Y\mid \mathcal{G}] = X$. 
    \end{enumerate}
\end{theorem}

\begin{proof}
    The proof is provided in Appendix~\ref{proof:umbrella-multi}.
\end{proof}

A direct observation from Theorem~\ref{th:umbrella-multi} is that the \ref{item:multi-cx} and \ref{item:multi-G} fail as soon as there is one $w$ not yielding an infinite mean in \eqref{eq:all-directions}. Hence, it is not sufficient for $\nu$'s marginals to all have undefined means; their dependence structure must be flexible enough to allow \eqref{eq:all-directions}. Independence, for instance, will always satisfy this flexibility requirement.





\subsection*{Acknowledgments}
 Ruodu Wang is supported by the Natural Sciences and Engineering Research Council of Canada (RGPIN-2018-03823, RGPAS-2018-522590) and Canada Research Chairs (CRC-2022-00141).

\appendix

\section{Proofs}
This appendix contains longer proofs which were omitted in the main text, namely the proofs of Lemma~\ref{lem:decompo} and Theorem~\ref{th:umbrella-multi}. Both proofs require Lemma~\ref{lem:ecx-multi}, which follows.

\subsection{Proof of Lemma~\ref{lem:decompo}}
\label{proof:decompo}

For any measure $\theta\in\mathcal{M}$, let $\pi_{\theta}$ be defined by $\pi_{\theta}(k) = \int (k-x)^+ \theta(\mathrm{d}x)$ and $c_{\theta}(k) = \int (x-k)^+ \theta(\mathrm{d}x)$ for every $k\in\mathbb{R}$. For measures $\theta_1,\theta_2$ with same total mass, the inequalities $\pi_{\theta_1}\leq \pi_{\theta_2}$ and $c_{\theta_1} \leq c_{\theta_2}$ characterize $\theta_1 \leq_{\rm cx} \theta_2$, even if $\theta_1,\theta_2$ do not have finite means; see \citet[Theorem~1]{CW26}.  The first part of the lemma's claim consists of the existence of a $\beta$ such that $\alpha \leq_{\rm cx} \beta \leq \eta$, and thus we rely on existing results on the existence of a shadow measure of $\alpha$ in $\eta$. However, such results from \cite{BJ16} and \cite{BHN22} assume that all considered distributions have finite means; the main difficulty of this proof is to argue that $\beta$ exists even without the  assumption of finite means. 

We argue three different cases: $\eta$ has a finite mean, a positive infinite mean, or an undefined mean. The case where it has a negative infinite mean follows by symmetry to the one with positive infinite mean.  

\textbf{For the case where $\eta$ has finite mean:} The assumptions $\lambda \leq_{\rm cx} \eta$ and $\alpha \leq \lambda$ directly imply that $\alpha \leq_{\rm ecx} \eta$. This is precisely the condition for the existence of a shadow measure \citep[Lemma~4.6]{BJ16}; let $\beta$ be that shadow measure, meaning $\alpha\leq_{\rm cx} \beta \leq \eta$. By \citet[Theorem~4.7]{BHN22}, 
\begin{equation*}
    \pi_{\beta} = \pi_{\eta} - (\pi_{\eta} - \pi_{\alpha})^{\rm c}
\end{equation*}
where $f^c$ means the greatest convex function below $f$. Thus, $\pi_{\eta - \beta} = \pi_{\eta} -\pi_{\beta} = (\pi_{\eta} - \pi_{\alpha})^{\rm c}$. 
Therefore, because $\lambda\leq_{\rm cx}\eta$ implies $\pi_{\lambda}\leq \pi_{\eta}$, it follows that $\pi_{\lambda - \alpha} = \pi_{\lambda} - \pi_{\alpha} \leq (\pi_{\eta} - \pi_{\alpha})^{\rm c}$ since $\pi_{\lambda-\alpha}$ is convex and below $\pi_{\eta} - \pi_{\alpha}$. Hence, $\pi_{\lambda - \alpha}\leq \pi_{\eta-\beta}$, 
and this yields $\lambda - \alpha \leq_{\rm cx} \eta - \beta$. 

\textbf{For the case where $\eta$ has positive infinite mean:} 
 First, we argue that there exists $m$ sufficiently large so that $\alpha\leq_{\rm ecx}\eta\vert_{(-\infty,m]}$. 
Take $q,m\in\mathbb{R}$, $m>q$, so that $\mathrm{supp}(\alpha)\subseteq(-\infty, q]$ and
\begin{equation*}
    \int(x-q)\eta^*(\mathrm{d}x) \geq \int(x-q)\alpha(\mathrm{d}x), 
\end{equation*}
where $\eta^* = \eta\vert_{(-\infty,m]}$ and $\ell = \eta^*(\mathbb{R}) - \alpha(\mathbb{R})$. Such $q$ and $m$ exist because $\alpha(\mathbb{R})<\eta(\mathbb{R})$ and $v \mapsto \int_{(-\infty,v]} x \eta(\mathrm{d}x)$ goes to infinity by the positive infinite mean assumption.  
Let $z = (\int x \eta^*(\mathrm{d}x) - \int x\alpha(\mathrm{d}x))/\ell$. We add mass at $z$ to $\alpha$
so that it has the same total mass as $\eta^*$; i.e., define $\alpha^*\in\mathcal{M}$ by letting $\alpha^* = \alpha + \ell \delta_{z}$. We note that $q\leq z \leq m$; indeed, $\pi_{\alpha}(m) \leq \pi_{\lambda}(m) \leq \pi_{\eta}(m) = \pi_{\eta^*}(m).$ Moreover, 
\begin{equation*}
    \pi_{\alpha^*}(k) = \pi_{\alpha}(k) \leq \pi_{\lambda}(k) \leq \pi_{\eta}(k) = \pi_{\eta^*}(k), \quad \text{for every }k\leq z;
\end{equation*}
and because $\alpha^*$ and $\eta^*$ have same total mass and means, 
\begin{equation*}
    \pi_{\alpha^*}(k) = k\eta^*(\mathbb{R}) - \int x\eta^*(\mathrm{d}x)  \leq \pi_{\eta^*}(k), \quad \text{for every }k\geq z.
\end{equation*}
Therefore, $\alpha^* \leq \eta^*$, because they have the same finite mean and $\pi_{\alpha^*}\leq \pi_{\eta^*}$; see \citet[Equation (3.A.6) and Theorem 3.A.1]{SS07}. Since $\alpha\leq \alpha^*$, we obtain $\alpha\leq_{\rm ecx} \eta^*$.

The measure $\eta^*$ has finite first moment, so this yields the existence of the shadow measure of $\alpha$ in $\eta^*$ by \citet[Lemma~4.6]{BJ16}.
Define $\beta$ as this shadow measure, so $\alpha\leq_{\rm cx} \beta \leq \eta^*$. This implies $\mathrm{supp}(\beta) \subseteq (-\infty, m]$, and hence $\beta$ has finite first moment. It remains to prove that $\lambda-\alpha\leq_{\rm cx} \eta - \beta$. Let $a>m$, and define $\lambda^{\wedge a}$ by moving all the mass beyond $a$ to $a$; i.e.,~if $X\sim \lambda$, then $\lambda^{\wedge a}$ is the distribution of $(X\wedge a)$. Similarly, let $b\geq a$, and define $\eta^{\wedge b}$ by moving all the mass beyond $b$ to $b$, but choose $b$ so that $\lambda^{\wedge a}$ and $\eta^{\wedge b}$ have the same mean. This is feasible because $b\mapsto \int(x\wedge b)\eta(\mathrm{d}x)$ is continuous and tends to infinity, and 
\begin{equation*}
\int(x\wedge a) \lambda(\mathrm{d}x) = a\lambda(\mathbb{R}) - \pi_{\lambda}(a) \geq a\eta(\mathbb{R}) - \pi_{\eta}(a) = \int(x\wedge a)\eta(\mathrm{d}x), 
\end{equation*}
where the inequality follows from $\pi_{\lambda}(a)\leq \pi_{\eta}(a)$ as $\lambda\leq_{\rm cx} \eta$. It is straightforward to see that $\pi_{\lambda^{\wedge a}} \leq \pi_{\eta^{\wedge b}}$.  Therefore, $\lambda^{\wedge a} \leq_{\rm cx} \eta^{\wedge b}$ as they have the same mean. 
Next, 
define $\beta^{\dagger}$ as the shadow measure of $\alpha$ in $\eta^{\wedge b}$; it exists because $\alpha\leq_{\rm ecx} \eta^{\wedge b}$ since $\eta* \leq \eta^{\wedge b}$ from $b>m$. 
We have $\alpha \leq \lambda^{\wedge a}$ and $\lambda^{\wedge a} \leq_{\rm cx}  \eta^{\wedge b}$, so we can apply the result already proved for the finite-mean case: $\beta^{\dagger}$ will satisfy both
\begin{equation}
    \alpha \leq_{\rm cx} \beta^{\dagger} \quad \text{and} \quad \lambda^{\wedge a} - \alpha \leq_{\rm cx} \eta^{\wedge b} - \beta^{\dagger}.\label{eq:beta-dagger}
\end{equation}
Next, we argue that $\beta = \beta^{\dagger}$. Note that $\beta \leq \eta^{\wedge b}$, so $\beta^{\dagger} \leq_{\rm cx} \beta$ by minimality of shadow measures with respect to the convex order. This means that $c_{\beta^{\dagger}}(m) \leq c_{\beta}(m) = 0$. Therefore, $\beta^{\dagger}$ has no mass above $m$.  Hence, $\beta^{\dagger} \leq \eta^*$, so $\beta \leq_{\rm cx} \beta^{\dagger}$ by minimality of shadow measures with respect to the convex order. This yields $\beta^{\dagger} = \beta$ because $\leq_{\rm cx}$ is antisymmetric on $L^1$. 
Finally, note that $\pi_{\lambda^{\wedge a}}(k) = \pi_{\lambda}(k)$  and $\pi_{\eta^{\wedge b}}(k) = \pi_{\eta}(k)$ for every $k\leq m$, so from \eqref{eq:beta-dagger}, $\pi_{\lambda-\alpha} (k) \leq\pi_{\eta - \beta}(k)$ for $k\leq m$. Also, when $k>m$,  $\pi_{\alpha}(k)=\pi_{\beta}(k)$, so $\pi_{\lambda-\alpha} (k) \leq\pi_{\eta - \beta}(k)$. As $\eta-\beta$ has infinite positive mean, hence $c_{\lambda-\alpha} (k) \leq c_{\eta - \beta}(k)$ holds for all $k\in\R$.
Therefore, $\lambda - \alpha \leq_{\rm cx} \eta - \beta$.

\textbf{For the case where $\eta$ has undefined mean:} 
By Lemma~\ref{lem:ecx-multi}, there exists $m_1,m_2\in\mathbb{R}$ such that $\alpha\leq_{\rm ecx} \eta\vert_{[m_1,m_2]}$. 
Since $\eta\vert_{[m_1,m_2]}$ has finite mean, we can invoke \citet[Lemma 4.6]{BJ16} for the existence of the shadow measure of $\alpha$ in $\eta\vert_{[m_1,m_2]}$. Let $\beta$ be that shadow measure, which has finite mean. Note that $\eta - \beta$ has undefined mean so $\lambda - \alpha \leq_{\rm cx} \eta - \beta$ holds by Lemma~\ref{lem:always}.

\subsection{Proof of Theorem~\ref{th:umbrella-multi}}
\label{proof:umbrella-multi}

This proof relies on Lemma~\ref{lem:ecx-multi}, whose statement and proof follow below (in Appendix~\ref{proof:ecx-multi}). 

For the implication \ref{item:multi-G} $\implies$ \ref{item:multi-cx}:  For any convex function $\varphi : \mathbb{R}^d\to \mathbb{R}$, 
\begin{equation*}
    \mathbb{E}[\varphi(X)] = \mathbb{E}\left[ \varphi\left( \mathbb{E}[Y\mid \mathcal{G}] \right)\right]  \leq \mathbb{E}\left[  \mathbb{E}[\varphi(Y)\mid \mathcal{G}]\right] = \mathbb{E}[\varphi(Y)],
\end{equation*}
provided both $\mathbb{E}[\varphi(X)]$ and $\mathbb{E}[\varphi(Y)]$ are well-defined, and where the inequality is Jensen's inequality. Hence, $X\leq_{\rm cx} Y$. 

For the implication \ref{item:multi-cx} $\implies$ \ref{item:multi-infinitemean}: Suppose that $\nu$ does not have infinite mean along every ray; this means that there exists $w_*\in\mathbb{R}^d\backslash\{0\}$ such that $\int_{\mathbb{R}^d} (w_*^{\mathrm{T}} y)_+ \nu(\mathrm{d}y) = k$ for some $k<\infty$. Choose $x\in\mathbb{R}^d$ such that $w_*^{\mathrm{T}} x> k$, and let $\mu$ be the Dirac measure at $x$. Then, 
$\int \varphi \,\mathrm{d}\mu > \int\varphi \,\mathrm{d}\nu$ for the convex function $y\mapsto (w_*^{\mathrm{T}} y)_+$, so $\mu\not\leq_{\rm cx} \nu$. 

For the implication \ref{item:multi-infinitemean} $\implies$  \ref{item:multi-G}:
Choose any $\mu\in\mathcal{P}_d$. Choose a compactly supported submeasure $\alpha_1\leq\mu$ such that $\alpha_1(\mathbb R^d)=\mu(\mathbb R^d)/2$.
By Lemma~\ref{lem:ecx-multi}, there exists $r>0$ such that $\alpha_1 \leq_{\rm ecx} \nu\vert_{B_r}$ where $B_r = \{y\in\mathbb{R}^d : ||y|| < r\}$, and write $\nu^* = \nu\vert_{B_r}$. By \citet[Theorem~4.2]{EH98}, $\alpha_{1}\leq_{\rm ecx} \nu^*$ implies the existence of $\gamma\in\mathcal{M}_d$ such that $\alpha_1\leq \gamma$ and $\gamma$ is a fusion of $\nu^*$; hence, $\gamma\leq_{\rm cx} \nu^*$. Both $\gamma$ and $\nu^*$ have finite means in each dimension, so Proposition~\ref{th:S1-multi} indicates that there are random vectors $U\sim \gamma/\gamma(\mathbb{R}^d)$ and $V\sim \nu^*/\nu^*(\mathbb{R}^d)$ such that $\mathbb{E}[V\mid \sigma(U)] = U$. Because $\alpha_1\leq \gamma$, the Radon-Nikodym derivative $\mathrm{d}\alpha_1/\mathrm{d}\gamma$  is well-defined and $0 \leq \mathrm{d}\alpha_1/\mathrm{d}\gamma \leq 1$. Define $\beta_1\in\mathcal{M}_d$ by letting 
\begin{equation*}
    \beta_1(A) = \gamma\left(\mathbb{R}^d\right)\mathbb{E}\left[\id_{\{V\in A\}} \frac{\mathrm{d}\alpha_1}{\mathrm{d}\gamma}(U) \right]\quad \text{for every }A\in\mathcal{B}\left(\mathbb{R}^d\right). 
\end{equation*}
Then, $\beta_1 \leq \nu^*$ and $\beta_1(\mathbb{R}^d) = \alpha_1(\mathbb{R}^d)$. Moreover, for every convex function $\varphi: \mathbb{R}^d \to \mathbb{R}$,  we have
\begin{align*}
    \int \varphi \mathrm{d}\alpha_1 &= \gamma\left(\mathbb{R}^d\right)\mathbb{E}\left[\varphi(U) \frac{\mathrm{d}\alpha_1}{\mathrm{d}\gamma}(U) \right]= \gamma\left(\mathbb{R}^d\right)\mathbb{E}\left[\varphi\left(\mathbb{E}[V\mid \sigma(U)] \right) \frac{\mathrm{d}\alpha_1}{\mathrm{d}\gamma}(U) \right]\\
&\leq\gamma\left(\mathbb{R}^d\right)\mathbb{E}\left[\mathbb{E}[\varphi(V)\mid \sigma(U)] \frac{\mathrm{d}\alpha_1}{\mathrm{d}\gamma}(U) \right]  = \gamma\left(\mathbb{R}^d\right)\mathbb{E}\left[\varphi(V)\frac{\mathrm{d}\alpha_1}{\mathrm{d}\gamma}(U) \right] = \int \varphi \mathrm{d}\beta_1, 
\end{align*}
where the inequality is Jensen's inequality and the last equality follows by extending the definition of $\beta_1$ from indicator functions to integrable functions using standard approximation. Therefore, $\alpha_1 \leq_{\rm cx} \beta_1$. Moreover, $\nu - \beta_1$ still has infinite mean in all directions since $\beta_1$ has a bounded support. Repeat this rationale to construct $\alpha_2$ and $\beta_2$ from $\mu - \alpha_1$ and $\nu - \beta_1$, and so on. Because half of the remaining mass is used at each iteration, $\sum_{n=1}^{\infty} \alpha_n(\mathbb{R}^d) = \mu(\mathbb{R}^d)$, and $\sum_{n=1}^{\infty} \beta_n(\mathbb{R}^d) = \nu(\mathbb{R}^d)$. Therefore 
\begin{equation*}
    \mu = \sum_{n=1}^{\infty} \alpha_n \quad \text{and} \quad\nu = \sum_{n=1}^{\infty} \beta_n
\end{equation*}
from the iterative construction. The argument in the proof of Theorem~\ref{th:S0} for the implication \ref{item:S0-alphabeta}$\implies$\ref{item:S0-G} then applies near identically, and yields the desired claim.

\subsection{A lemma on the multivariate extended convex order}
\label{proof:ecx-multi}
The multivariate analogue to $\leq_{\rm ecx}$ defined above is as follows: For $\mu,\nu\in\mathcal{M}_d$, we write $\mu\leq_{\rm ecx} \nu$ when $\mu(\mathbb{R^d})\leq \nu(\mathbb{R^d})$ and $\int \varphi \,\mathrm{d}\mu \leq \int \varphi \,\mathrm{d}\nu$ holds for every nonnegative convex function $\varphi:\mathbb{R}^d \to \mathbb{R}$. 

\begin{lemma}
\label{lem:ecx-multi}
    Consider $\alpha,\eta\in\mathcal{M}_d$ such that $\eta$ has infinite mean along every ray, $\alpha$ has compact support and $\alpha(\mathbb{R}^d) < \eta(\mathbb{R}^d)$. Then, there is a ball $B_r = \{y\in\mathbb{R}^d:||y||<r\}$ such that 
    \begin{equation*}
        \alpha \leq_{\rm ecx} \eta\vert_{B_r}. 
    \end{equation*}
\end{lemma}
\begin{proof}
    Let $m>0$ satisfy $\mathrm{supp}(\alpha)\subseteq B_m$. We suppose that $\alpha \leq_{\rm ecx} \eta\vert_{B_n}$ fails for every integer $n\geq m$ and show that this yields a contradiction. For every $n\geq m$, the failure of $\alpha\leq_{\rm ecx} \eta\vert_{B_n}$ means that we can pick nonnegative convex $\varphi_n$ satisfying $\int_{B_n}\varphi_n\,\mathrm{d}\eta < \int\varphi_n \,\mathrm{d}\alpha$. At a point maximizing $\varphi_n$ on $\bar{B}_m$, where $\bar{B}_m$ denotes the closure of $B_m$, let $\ell_n$ be an affine lower bound touching $\varphi_n$, meaning $\ell_n$ is in the form $\ell_n(y) = a^{\mathrm{T}}_n y + b_n$ and $\sup_{B_m}\ell_n = \sup_{B_m} \varphi_n>0$. The positivity of $\sup_{B_m} \varphi_n$ follows from $\int\varphi_n \,\mathrm{d}\alpha>0$. Therefore, 
    \begin{equation}
        \int_{B_n}(a_n^{\mathrm{T}} y + b_n)_+\eta(\mathrm{d}y)\leq\int_{B_n}\varphi_n\,\mathrm{d}\eta<\int\varphi_n\,\mathrm{d}\alpha \leq\alpha(\R^d)\sup_{y\in B_m} \ell_n(y)\leq \alpha(\mathbb{R}^d)(m||a_n|| + b_n).
        \label{eq:ecx-1}
    \end{equation}
    For each $n\geq m$, note $\sup_{B_m}\ell_n>0$, so $\sqrt{||a_n||^2+b_n^2}>0$. Let $a_n^*, b_n^*$ be a rescaling of $a_n,b_n$ so that $||a_n^*||^2 + (b_n^*)^2 = 1$. There exists a subsequence $\{(a_{n_j}^*, b_{n_j}^*):j\in\mathbb{N}\}$ where $(a^*_{n_j}, b^*_{n_j}) \to (a,b)$ as $j\to\infty$ for some $a\in\mathbb{R}^d$, $b\in\mathbb{R}$. Note that $\lim_{j\to\infty}\id_{B_{n_j}}(y)(a_{n_j}^{*\mathrm{T}}y + b_{n_j}^*)_+=(a^{\mathrm{T}}y + b)_+$. By \eqref{eq:ecx-1} and Fatou's lemma, 
    \begin{equation}
        \int(a^{\mathrm{T}}y + b)_+ \eta(\mathrm{d}y) \leq \alpha(\mathbb{R}^d)(m||a|| + b) < \infty.
        \label{eq:ecx-2}
    \end{equation}
 Since $\eta$ has infinite mean along every ray, we must have $a=(0,\ldots,0)$, and hence $b\neq0$. If $b>0$, \eqref{eq:ecx-2} implies $\eta(\mathbb R^d)\leq\alpha(\mathbb R^d)$, contradicting our assumption. If $b<0$, it implies $\alpha(\mathbb R^d)=0$, contradicting $\int\varphi_n \,\mathrm{d}\alpha>0$. Thus, both cases lead to a contradiction.
\end{proof}



\end{document}